# Dutch book in simple multivariate normal prediction: Another look

## Morris L. Eaton[1]

*University of Minnesota*

**Abstract:** In this expository paper we describe a relatively elementary method of establishing the existence of a Dutch book in a simple multivariate normal prediction setting. The method involves deriving a nonstandard predictive distribution that is motivated by invariance. This predictive distribution satisfies an interesting identity which in turn yields an elementary demonstration of the existence of a Dutch book for a variety of possible predictive distributions.

## 1. Introduction

Ordinarily, showing that a popular inferential scheme suffers from de Finetti's incoherence (existence of a Dutch book) is met with surprise since incoherence is typically regarded as a serious indictment of a statistical method. An instance of this incoherence occurs in a simple $p$-dimensional multivariate normal setting. The statistical problem is to predict the next observation in a sequence of independent and identically distributed normal observations (with mean 0 and unknown $p \times p$ positive definite covariance matrix $\Sigma$). The "usual" predictive distribution is obtained from a formal Bayes calculation using the Jeffreys' invariant prior distribution for $\Sigma$. This is detailed in Eaton and Sudderth [9] where it is shown that the "usual" solution is in fact incoherent (a Dutch book can be made against the "usual" predictive distribution). The arguments in Eaton and Sudderth [9] are neither simple nor intuitive, but variations on these have been used effectively to extend the Dutch book argument to other multivariate settings (see Eaton and Sudderth [10, 11, 12, 13]). Some background information on incoherence and Dutch book is given in Section 2 below.

Recently, Eaton and Freedman [8] presented a relatively simple and self-contained argument that the Jeffreys' invariant prior resulted in a Dutch book in the simple normal prediction problems. The demonstration relies mainly on sampling properties of the normal distribution rather than on invariance arguments (see Eaton and Sudderth [9]). However, the Eaton and Freedman [8] arguments seem to be rather special and not easily adaptable to other invariant proposals such as those discussed in Bjørnstad [2].

The focus of this paper is again the simple normal prediction problem of Eaton and Sudderth [9] and Eaton and Freedman [8]. The purpose of writing the current paper is to present an argument that:

[1]School of Statistics, University of Minnesota, 367 Ford Hall, 224 Church Street S.E., Minneapolis, MN 55455, USA, e-mail: eaton@stat.umn.edu

(i) reveals the role of the invariance in the incoherence (existence of a Dutch book);

(ii) relies primarily on calculus so is mainly self-contained;

(iii) yields the Dutch book conclusion for many predictive proposals in the simple normal prediction problem.

To put things into perspective it is useful to sketch the argument used below. As in Eaton and Sudderth [9], let $G_T^+$ be the group of $p \times p$ lower triangular matrices with positive diagonal elements. Then the unknown covariance matrix $\Sigma$ in the normal model can be uniquely expressed as $\Sigma = \theta\theta'$, $\theta \in G_T^+$. Using the right Haar measure on $\theta \in G_T^+$ as an improper prior distribution (rather than the Jeffreys' prior on $\Sigma$), a formal Bayes calculation yields a predictive distribution that we denote by $Q_H$. The predictive distribution $Q_H$ is invariant under the group $G_T^+$. Details and a full explanation are given in Section 3.

The remainder of the paper is devoted to an argument that allows one to compare almost any $G_T^+$ invariant predictive distribution $Q$ with the special predictive distribution $Q_H$ above. The argument is based on recent work of Zhu [21] and on a density based approach of Eaton and Sudderth [13]. In essence, the constructive method described in Sections 4 and 5 shows that if an invariant predictive distribution $Q$ has a density $q$, and if $q$ differs on a set of positive Lebesgue measure from the density of $Q_H$, then $Q$ is incoherent (a Dutch book exists). This result is used to show a Dutch book exists for several well-known predictive proposals.

Finally, we should mention that $Q_H$ is not incoherent because the group $G_T^+$ is amenable. See Eaton and Sudderth [11], Theorem 8.1, for a general discussion and a proof.

## 2. Background

The focus of this paper concerns a method for the evaluation of predictive distributions in a simple multivariate normal sampling situation. However, it is useful to first describe what we mean by a prediction problem and to detail the evaluative criterion of interest here.

The origin of our formulation of the prediction problem stems from Laplace [17] (see Stigler [19] for an English translation of the original French). Consider a sample space $(\mathcal{X}, B_1)$ and an observation (usually a vector or a matrix) $X$ in $\mathcal{X}$. A variable $Z$ taking values in $\mathcal{Z}$ is to be predicted on the basis of an assumed joint parametric probability model

$$P(dx, dz|\theta), \quad \theta \in \Theta.$$

Here $\Theta$ is a parameter space and $\theta$ is an unknown parameter. By a <u>predictive distribution</u> we mean a probability distribution $Q(dz|x)$ for $Z$ that is allowed to depend on the observed value of $X = x$. The primary example of concern throughout this paper is the following.

**Example 1.** Let $X_1, \ldots, X_n$ be independent and identically distributed (iid) with a $p$-dimensional multivariate normal distribution $N_p(0, \Sigma)$ with mean 0 and unknown $p \times p$ positive definite covariance matrix $\Sigma$. It is assumed that $n \geq p$, and $p \geq 2$. The sample space $\mathcal{X}$ is the set of all $p \times n$ matrices with rank $p$ (a set of Lebesgue measure zero has been discarded, for convenience). Thus, the matrix

$$X = (X_1, \ldots, X_n) : \ p \times n$$

is the "data point" in $\mathcal{X}$ and

$$(2.1) \qquad S = XX' = \sum_{i=1}^{p} X_i X_i' : \quad p \times p$$

has rank $p$. The variable $Z$ to be predicted is assumed to be $N_p(0, \Sigma)$ and independent of the $X$'s. One naive predictive distribution for $Z$ is

$$(2.2) \qquad Q_1(\cdot | x) = N_p(0, n^{-1}s),$$

where $x$ is the observed value of $X$ and $s$ is the observed value of $S$. Of course, the intuition is that $n^{-1}S$ is an unbiased estimator of $\Sigma$ based on a minimal sufficient statistic. An alternative proposal, based on a formal Bayes argument and the so-called Jeffreys' improper prior distribution (see Eaton and Freedman [8] for some discussion and details), yields

$$(2.3) \qquad Q_2(dz | x) = q_2(z | x) dz,$$

where $dz$ is Lebesgue measure on $\mathcal{Z} = R^p$ and the density $q_2(\cdot | x)$ is

$$(2.4) \qquad q_2(z | x) = C_{n,p} \frac{|s|^{-\frac{1}{2}}}{(1 + z's^{-1}z)^{(n+1)/2}}.$$

The constant $C_{n,p}$ is

$$(2.5) \qquad C_{n,p} = \frac{\Gamma(\frac{n+1}{2})}{\pi^{p/2} \Gamma(\frac{n-p+1}{2})}$$

and $|s|$ is the determinant of $s$. This ends the introduction of Example 1.

We now return to the general prediction setting at the beginning of this section. Stone [20] introduced the notion of <u>strong</u> <u>inconsistency</u> as a criteria for excluding certain types of probability distributions in inferential settings. In the prediction framework, this idea takes the following form.

**Definition 2.1.** A <u>predictive distribution</u> $Q(dz | x)$ is <u>strongly inconsistent</u> (SI) with the model $\{P(dx, dz | \theta) | \theta \in \Theta\}$ if there exists a measurable function $\bar{f}(x, z)$ with values in $[-1, 1]$ and an $\epsilon > 0$ so that

$$(2.6) \qquad \sup_x \int f(x, z) Q(dz | x) + \epsilon \leq \inf_\theta \int \int f(x, z) P(dx, dz | \theta).$$

The intuition behind SI is that when (2.6) holds, no matter what the distribution for $X$, say $m(dx)$,

$$\int \int f(x, z) Q(dz | x) m(dx) + \epsilon \leq \int \int f(x, z) P(dx, dz | \theta),$$

for all $\theta$. Thus, if $Q(dz | x)$ is used as a distribution for $Z$ after seeing $X = x$, then under all possible models for $(X, Z)$ that are consistent with $Q(dz | x)$, the expectation of $f$ is at least $\epsilon$ less than any expectation of $f$ under the assumed model. Hence the terminology, strong inconsistency.

Ramsay [18] and independently de Finetti [3, 4] introduced the notion of betting schemes for the evaluation of proposed probability distributions. Their ideas were extended by Freedman and Purves [14] to cover cases with conditional bets. A somewhat modified approach, due to Heath and Sudderth [15], is relevant to the discussion below. In the prediction case when the spaces are infinite (as in Example 1), the Ramsay–de Finetti–Freedman and Purves scheme takes the following form. Let $C$ be a subset of $\mathcal{X} \times \mathcal{Z}$ and let $C_x = \{z|(x, z) \in C\}$ be the $x$-section of $C$. If an inferrer is using $Q(dz|x)$ as a predictive distribution and $X = x$ is observed, then $C_x \subseteq \mathcal{Z}$ is assigned the probability $Q(C_x|x)$. Therefore

$$\Psi_0(x, z) = I_C(x, z) - Q(C_x|x)$$

has $Q(\cdot|x)$ – expectation zero. A standard interpretation of $\Psi_0$ as a payoff function is:

> After seeing $X = x$, the inferrer gives $C_x$ probability $Q(C_x|x)$. A gambler can pay $Q(C_x|x)$ - dollars for a ticket worth
>
> $$\left\{ \begin{array}{ll} \$1 & \text{if } Z \in C_x \\ 0 & \text{if } Z \notin C_x. \end{array} \right.$$
>
> Obviously, the net payoff to the gambler is $\Psi_0(x, z)$. The inferrer regards the bet as "fair" since $\Psi_0$ has expectation zero under $Q(\cdot|x)$.

Slightly more complicated betting scenarios are constructed as follows. Suppose the gambler can pick subsets $C_1, \ldots, C_r$ in $\mathcal{X} \times \mathcal{Z}$ and pays $c_i(x)Q(C_{i,x}|x)$ dollars for a ticket worth

$$\left\{ \begin{array}{ll} \$c_i(x) & \text{if } Z \in C_{i,x} \\ 0 & \text{if } Z \notin C_{i,x}, \end{array} \right.$$

for $i = 1, \ldots, r$. The numbers $c_i(x)$ are assumed to be bounded functions of $x$, but need not be non-negative. The net payoff to the gambler is computed by summing the individual payoffs, so the <u>net payoff function</u> is

$$(2.7) \qquad \Psi(x, z) = \sum_{i=1}^{r} c_i(x) \left[ I_{C_i}(x, z) - Q(C_{i,x}|x) \right].$$

Again, since $\Psi(x, \cdot)$ has $Q(\cdot|x)$ expectation 0, the inferrer regards the betting scheme of the gambler as fair.

Here is de Finetti's notion of incoherence adapted to the current prediction setting.

**Definition 2.2.** The predictive distribution $Q(\cdot|x)$ is <u>incoherent</u> if there is an $\epsilon > 0$ and a payoff function $\Psi$ of the form (2.7) such that

$$(2.8) \qquad \mathcal{E}_\theta \Psi(X, Z) \geq \epsilon \text{ for all } \theta \in \Theta.$$

In other words, the predictive distribution is incoherent if the gambler has a uniformly (over $\theta$) positive expected gain under the model. The discussion leading to Definition 2.2 comes from Heath and Sudderth [15]. The Freedman and Purves [14] formulation was in terms of odds rather than payoff functions.

In examples, the verification of SI or incoherence is typically not straightforward. However, as shown recently in Eaton and Freedman [8], in the prediction setting of this paper, the two notions are equivalent (see Theorem 2 on p. 868-869). For the example under consideration here, SI is established for a variety of predictive distributions (see Section 5).

Finally, the term "Dutch book" is used as a synonym for incoherence. When incoherence obtains, standard terminology is to say "Dutch book" can be made against the predictive distribution $Q(\cdot|x)$. See Eaton and Freedman [8] for some discussion and a bit of history.

### 3. The Haar predictive distribution

As in Example 1, consider $X_1, \ldots, X_n$ iid $N_p(0, \Sigma)$ where $\Sigma$ is an unknown $p \times p$ positive definite matrix. It is assumed that $n \geq p$ so the matrix $S$ in (2.1) is positive definite when the sample matrix

$$X = (X_1, \ldots, X_n): \ p \times n$$

is an element of $\mathcal{X}$. A variable $Z \in R^p \equiv \mathcal{Z}$ which is $N_p(0, \Sigma)$ is to be predicted after seeing $X = x \in \mathcal{X}$. Here $Z$ is independent of $X_1, \ldots, X_n$. In short, the statistical problem of concern in this paper is to produce a predictive distribution $Q(dz|x)$ for $Z \in R^p$ after seeing the data $X = x$. The focus of this section is on a particular predictive distribution obtained via a formal Bayes calculation.

Recall that $G_T^+$ is the group of $p \times p$ lower triangular matrices with positive diagonal elements. For use in what follows, we list some facts about $G_T^+$ (see Eaton [7], especially Chapter 1 and in particular, pages 18 and 19). Elements $g \in G_T^+$ have positive diagonals, $g_{ii}, i = 1, \ldots, p$ and $g_{ij} = 0$ for $i < j$. The symbol "$dg$" denotes Lebesgue measure on $G_T^+$ (as an obvious open subset of $R^{p(p+1)/2}$). The measure

$$(3.1) \qquad \nu_r(dg) = \frac{dg}{\prod\limits_{i=1}^{p} g_{ii}^{p-i+1}}$$

is a <u>right invariant</u> (Haar) <u>measure</u> on $G_T^+$. The function

$$(3.2) \qquad \Delta(g) = \prod_{i=1}^{p} g_{ii}^{p-2i+1}$$

is the <u>modular function</u> of $G_T^+$ and

$$(3.3) \qquad \nu_l(dg) = \Delta(g)\nu_r(dg) = \frac{dg}{\prod\limits_{i=1}^{p} g_{ii}^{i}}$$

is a <u>left invariant measure</u> on $G_T^+$.

Given a $p \times p$ positive definite matrix $E$, there is a unique element $T \in G_T^+$ such that $E = TT'$ (see Eaton [6], Proposition 5.4 for a proof). This element is denoted by $\tau(E)$ in some expressions below. In particular, $\theta = \tau(\Sigma)$ is a reparameterization of covariance matrices. In this parameterization, the density function of $X$, with respect to Lebesgue measure on $\mathcal{X}$, is

$$(3.4) \qquad f_1(x|\theta) = \frac{|\theta|^{-n}}{(2\pi)^{np/2}} \exp\left\{-\frac{1}{2} tr(\theta\theta')^{-1}s\right\}, \ \ x \in \mathcal{X}$$

where $s = xx': p \times p$, $\theta \in G_T^+$, and "tr" denotes trace. Of course the density function of $Z$ on $R^p$ is

$$(3.5) \qquad f_2(z|\theta) = \frac{|\theta|^{-1}}{(2\pi)^{p/2}} \exp\left\{-\frac{1}{2} tr(\theta\theta')^{-1}zz'\right\}.$$

As usual, $|\cdot|$ denotes the determinant. These two densities define the probability models

$$(3.6) \qquad P_1(dx|\theta), \ \ \theta \in \Theta = G_T^+, \ \ x \in \mathcal{X}$$

and

$$(3.7) \qquad P_2(dz|\theta), \ \ \theta \in \Theta, \ \ z \in R^p,$$

for $X$ and $Z$ respectively. Recall $X$ and $Z$ are assumed independent so the joint model for $(X, Z)$ is

$$(3.8) \qquad P(dx, dz|\theta) = P_1(dx|\theta)P_2(dz|\theta)$$

and the joint density is

$$(3.9) \qquad f(x, z|\theta) = f_1(x|\theta)f_2(z|\theta), \ \ (x, z) \in \mathcal{X} \times \mathcal{Z}.$$

For $g \in G_T^+$, it is easy to show that

$$(3.10) \qquad f_1(gx|g\theta) = |g|^{-n}f_1(x|\theta)$$

and

$$(3.11) \qquad f_2(gz|g\theta) = |g|^{-1}f_2(x|\theta).$$

These two invariance properties of the densities imply that for the model (3.8),

$$(3.12) \qquad P(gB|g\theta) = P(B|\theta),$$

for all Borel sets $B \subseteq \mathcal{X} \times \mathcal{Z}$ and all $\theta$, $g \in G_T^+$. In other words, the assumed statistical model is $G_T^+$ invariant. Note that the group $G_T^+$ acts transitively on $\Theta = G_T^+$. In such invariant settings, it was argued in Eaton and Sudderth [13] that the use of the right Haar measure as a prior distribution may yield predictive distributions with interesting inferential properties.

We now proceed to calculate, via a formal use of Bayes Theorem, the predictive distribution (henceforth called the Haar inference) induced by using $\nu_r(d\theta)$ in (3.1) as an improper prior. To this end, let

$$(3.13) \qquad m(x, z) = \int_{G_T^+} f(x, z|\theta)\nu_r(d\theta)$$

and

$$(3.14) \qquad m_1(x) = \int_{G_T^+} f_1(x|\theta)\nu_r(d\theta),$$

for $x \in \mathcal{X}$ and $z \in R^p$. For fixed $x$,

$$(3.15) \qquad q_H(z|x) = m(x, z)/m_1(x)$$

is a density on $R^p$ and by definition, the predictive distribution $Q_H(dz|x)$ with density (3.15) is the <u>Haar inference</u>. Using (3.10), (3.11), and the properties of $\nu_r$, it is easy to show that $q_H$ is invariant in the sense that

$$(3.16) \qquad q_H(gz|gx) = |g|^{-1}q_H(z|x).$$

From (3.16), the invariance of $Q_H$,

$$(3.17) \qquad Q_H(gB|gx) = Q_H(B|x)$$

is immediate.

Because the calculation is not quite standard, we sketch the details in the derivation of $q_H$ in (3.15). For $w \in R^p$, define the function $\psi_p(w)$ by

$$(3.18) \qquad \psi_p(w) = \begin{cases} 1 & \text{if } p = 1 \\ (1 + w'w)^{-(p-1)/2} \prod_{i=1}^{p-1} (1 + w_1^2 + \ldots + w_i^2) & \text{if } p \geq 2. \end{cases}$$

**Lemma 3.1.** *For $w \in R^p$, recall $\tau(I_p + ww')$ is the unique element in $G_T^+$ that satisfies $(I_p + ww') = [\tau(I_p + ww')][\tau(I_p + ww')]'$. Then,*

$$(3.19) \qquad \Delta(\tau(I_p + ww')) = \psi_p(w)$$

*where $\Delta$ is the modular function given in (3.2).*

*Proof.* The proof of this (a messy calculation) is not too hard via an induction argument on dimension $p$. The details are omitted. $\qquad \square$

**Theorem 3.1.** *Let $k_0$ be the density on $R^p$ given by (see 2.4)*

$$(3.20) \qquad k_0(w) = \frac{\Gamma(\frac{n+1}{2})}{\pi^{p/2}\Gamma(\frac{n-p+1}{2})} \; \frac{1}{(1 + w'w)^{(n+1)/2}}.$$

*Then*

$$(3.21) \qquad k_1(w) = k_0(w)\frac{1}{\psi_p(w)}$$

*is a density on $R^p$ and*

$$(3.22) \qquad q_H(z|x) = |L|^{-1}k_1(L^{-1}z)$$

*where $xx' = s = LL'$ with $L \in G_T^+$.*

*Proof.* Using the expressions (3.4) and (3.5) and the fact that $\nu_r$ transforms to $\nu_l$ under the mapping $\theta \longrightarrow \theta^{-1}$ in $G_T^+$, it follows from (3.15) that

$$(3.23) \qquad q_H(z|x) = (2\pi)^{-p/2}\frac{\int |\theta|^{n+1}\exp\{-\frac{1}{2}tr\theta'\theta[s + zz']\}\nu_l(d\theta)}{\int |\theta|^n \exp\{-\frac{1}{2}tr\theta'\theta s\}\nu_l(d\theta)}.$$

Writing $s = LL'$ with $L \in G_T^+$ and setting $w = L^{-1}z$, some algebra and a change of variable yields

$$q_H(z|x) = (2\pi)^{-p/2}\frac{\int |\theta|^{n+1}\exp\{-\frac{1}{2}tr(\theta L)'(\theta L)(I_p + ww')\}\nu_l(d\theta)}{\int |\theta|^n \exp\{-\frac{1}{2}tr(\theta L)'(\theta L)\}\nu_l(d\theta)}$$

$$= |L|^{-1}(2\pi)^{-p/2}\frac{\int |\theta|^{n+1}\exp\{-\frac{1}{2}tr\theta'\theta(I_p + ww')\}\nu_l(d\theta)}{\int |\theta|^n \exp\{-\frac{1}{2}tr\theta'\theta\}\nu_l(d\theta)}.$$

Setting $I_p + ww' = UU'$, $U \in G_T^+$ and changing variables in the numerator integral above gives

$$(3.24) \qquad q(z|x) = \frac{|L|^{-1}|U|^{-(n+1)}}{(2\pi)^{p/2}\Delta(U)} \cdot \frac{\int |\theta|^{n+1}\exp\{-\frac{1}{2}tr\theta'\theta\}\nu_l(d\theta)}{\int |\theta|^n \exp\{-\frac{1}{2}tr\theta'\theta\}\nu_l(d\theta)}.$$

Now, note that

$$|U|^{-(n+1)} = |UU'|^{-\frac{n+1}{2}} = \frac{1}{(1+w'w)^{(n+1)/2}}$$

and from Lemma 3.1,

$$\Delta(U) = \psi_p(w).$$

Finally, a standard (but not routine) multivariate calculation yields

$$\frac{1}{(2\pi)^{p/2}} \frac{\int |\theta|^{n+1} \exp\{-\frac{1}{2}tr|\theta'\theta\} \nu_l(d\theta)}{\int |\theta|^n \exp\{-\frac{1}{2}tr|\theta'\theta\} \nu_l(d\theta)} = \frac{\Gamma(\frac{n+1}{2})}{\pi^{p/2}\Gamma(\frac{n-(p+1)}{2})}.$$

Piecing all of this together yields the expression (3.22) for $q_H(z|x)$. The fact that $k_1$ is a density on $R^p$ follows by setting $L = I_p$ in (3.22). □

Let $q(z|x) \geq 0$ be an arbitrary predictive density for $Z$ with $X = x$. That is, $q(\cdot|x)$ is a density on $R^p$ for each $x \in \mathcal{X}$ and $q$ is jointly measurable on $\mathcal{Z} \times \mathcal{X}$.

**Definition 3.1.** The density $q(z|x)$ is $G_T^+$-invariant if for all $g \in G_T^+$,

$$(3.25) \qquad q(gz|gx) = |g|^{-1}q(z|x) \quad \text{for all } z, x.$$

Each $G_T^+$-invariant predictive density yields a $G_T^+$-invariant predictive distribution $Q(dz|x)$ given by

$$Q(B|x) = \int_B q(z|x)dz, \quad B \subseteq R^P.$$

The invariance of $Q$, namely $Q(gB|gx) = Q(B|x)$, is immediate from (3.25).

The density of $Q_1$ in (2.2) and $Q_2$ in (2.3) are both $G_T^+$-invariant. Further, when specialized to the case of mean zero considered in this paper, the first seven entries in Table 1 of the survey paper of Keyes and Levy [16] are all $G_T^+$-invariant. It is such predictive densities that are compared to $q_H$ in Section 5.

## 4. Zhu's result

Using some results of Eaton and Sudderth [12], Zhu [21] was able to establish an interesting and useful relationship between an invariant prediction model and the predictive distribution obtained from the right Haar measure. Zhu's result will be stated here only for the normal model under consideration. For a simplified version and proof of this general result when densities exist, see Eaton and Sudderth [13]. The most general version is Theorem 3.4.1 in Zhu [21].

In the notation of Section 3, let $P_1(dx|\theta)$, $P_2(dz|\theta)$ and $P(dx, dz|\theta)$ be the probability measures in (3.6), (3.7) and (3.8) respectively. Also, let $Q_H(dz|x)$ be the Haar inference defined by the density $q_H(z|x)$ in (3.15). Next, introduce the <u>Haar model</u> given by the probability measure

$$(4.1) \qquad P_H(dx, dz|\theta) = Q_H(dz|x)P_1(dx|\theta)$$

on $\mathcal{X} \times \mathcal{Z}$.

Recall that a real valued function $f$ on $\mathcal{X} \times \mathcal{Z}$ is <u>$G_T^+$-invariant</u> if for all $g \in G_T^+$

$$f(gx, gz) = f(x, z) \quad \text{for all } x, z.$$

**Theorem 4.1.** *The original model $P(dx, dz|\theta)$ and the Haar model $P_H(dx, dz|\theta)$ agree on the bounded invariant functions. That is, for each bounded measurable $G_T^+$-invariant $f$ and for all $\theta$,*

$$(4.2) \qquad \int \int f(x, z) P(dx, dz|\theta) = \int \int f(x, z) P_H(dx, dz|\theta).$$

The implications of (4.2) in more general settings are discussed in Eaton and Sudderth [13]. Because densities exist in the setting of this paper, Theorem 3.1 in Eaton and Sudderth [13] applies directly to this case so we omit a proof.

**Remark 1.** In the proof of Theorem 3.1 in Eaton and Sudderth [13], there are two misprints: On page 500 of this paper in lines 8 and 9 from the top, "$p(x, z|g\theta)$" should be "$p(x, z|\theta)$."

The use of (4.2) in this paper occurs in the next section which deals with the Dutch book argument.

## 5. Dutch book

Consider a predictive density $q(z|x)$ that is $G_T^+$-invariant. The purpose of this section is to show that if $q(z|x)$ is essentially different from $q_H(z|x)$, then the predictive distribution $Q(dz|x)$ determined by $q$ is incoherent so a Dutch book exists. The construction of the required pay-off function is explicit and both (2.6) and (2.8) are verified directly.

To make the above precise, consider the set

$$(5.1) \qquad C = \{(x, z)|q(z|x) < q_H(z|x)\} \subseteq \mathcal{X} \times \mathcal{Z}.$$

Then

$$(5.2) \qquad C_x = \{z|(x, z) \in C\}$$

is the $x$-section of the set $C$. Let $l$ denote Lebesgue measure on $\mathcal{X} \times \mathcal{Z}$, $l_1$ denote Lebesgue measure on $\mathcal{X}$ and $l_2$ denote Lebesgue measure on $\mathcal{Z}$. Clearly, $l(dx, dz) = l_1(dx)l_2(dz)$. From Tonelli's Theorem (see Dunford and Schwartz [5], p. 194),

$$(5.3) \qquad l(C) = \int_{\mathcal{X}} l_2(C_x) l_1(dx).$$

**Theorem 5.1.** *If $l(C) > 0$, then $Q(dz|x)$ is incoherent (strongly inconsistent) and a Dutch book exists for the predictive distribution $Q(\cdot|x)$.*

*Proof.* When $l(C) > 0$, the set

$$(5.4) \qquad D = \{x|l_2(C_x) > 0\} \subseteq \mathcal{X}$$

must have positive $l_1$-measure. Consider the function

$$(5.5) \qquad \phi(x, z) = I_D(x)[I_{C_x}(z) - Q(C_x|x)].$$

We first claim that the function $\phi(x, z)$ is $G_T^+$-invariant. To see this, first note that $C$ is an invariant set from the invariance of $q_H$ and the assumed invariance of $q$. From this it follows easily that

$$C_{gx} = gC_x, \quad g \in G_T^+.$$

Therefore the set $D$ is invariant and $x \longrightarrow Q(C_x|x)$ is an invariant function of $x$. The invariance of $\phi$ in (5.5) is now immediate.

We now verify (2.6) with $f = \phi$. First, for each $x$,

$$\int \phi(x, z) Q(dz|x) = 0.$$

Since $\phi$ is invariant and the model is invariant,

(5.6) $$\theta \longrightarrow \int \int \phi(x, z) P(dx, dz|\theta)$$

is an invariant function of $\theta$. Because $\Theta = G_T^+$, the right side of (5.5) is constant in $\theta$, say

$$\epsilon_0 = \int \int \phi(x, z) P(dx, dz|\theta) \quad \text{for all } \theta.$$

Now, apply Theorem 4.1 to obtain

$$\begin{aligned}
\epsilon_0 &= \int \int \phi(x, z) P(dx, dz|\theta) \\
&= \int \int \phi(x, z) P_H(dx, dz|\theta) \\
&= \int \int I_D(x)[I_{C_x}(z) - Q(C_x|x)] Q_H(dz|x) P_1(dx|\theta) \\
&= \int_D \int_{C_x} [q_H(z|x) - q(z|x)] l_2(dz) f_1(x|\theta) l_1(dx)
\end{aligned}$$

where $f_1$ is given in (3.4). This set $D$ has positive $l_1$-measure, and the density $f_1$ is positive everywhere. Also, for all $x \in D$, $C_x$ has positive $l_2$ measure and for all $z \in C_x$, $q_H(z|x) > q(z|x)$. Thus $\epsilon_0 > 0$ and (2.6) holds with $\epsilon = \epsilon_0$. Hence SI holds.

To see that incoherence holds, take the net payoff function to be $\phi(x, z)$ (in (2.7), take $r = 1$, $c_1(x) = I_D(x)$ and $C_1 = C$). A repeat of the argument above shows (2.8) holds with $\epsilon = \epsilon_0$. This completes the proof. □

Before discussing any examples, it is useful to next consider the case when the set $C$ has measure zero.

**Theorem 5.2.** *If $l(C) = 0$, then the set $D$ in (5.4) has $l_1$-measure zero and*

(5.7) $$Q(\cdot|x) = Q_H(\cdot|x) \quad a.e. \ (l_1).$$

*Proof.* When $l(C) = 0$, (5.3) shows $l_2(C_x) = 0$ a.e. $(l_1)$. Hence $D$ has $l_1$-measure 0. But for $x \in D^c$, $l_2(C_x) = 0$. For this $x$, a standard argument shows that $Q(\cdot|x) = Q_H(\cdot|x)$. Therefore (5.7) holds and the proof is complete. □

We now proceed to give a wide class of examples where the set $C$ has positive Lesbesgue measure so Theorem 5.1 applies. To this end, let $k$ be a density function on $R^p$ and as usual, for $x \in \mathcal{X}$ write

$$xx' = LL' \ , \ \ L \in G_T^+.$$

Observe that the predictive density $q_k(z|x)$ given by

(5.8) $$q_k(z|x) = |L|^{-1} k(L^{-1}z)$$

is $G_T^+$-invariant and is obviously determined by $k$. Note that the predictive distributions (2.2) and (2.3) both have predictive densities of the form (5.8) for an appropriate $k$. Further, the Haar predictive distribution $Q_H$ derived in Section 3 has the form (5.8) with $k = k_1$ where $k_1$ is given in (3.21).

The next result shows that any predictive distribution of the form (5.8) is incoherent when $k$ is different (on a set of positive $l_2$-measure) from the special density $k_1$ that defines the Haar inference $Q_H$.

**Theorem 5.3.** *Let $k$ be a density on $R^p$ and let $k_1$ be the density defined in (3.21). Assume that*

$$\int |k(z) - k_1(z)| dz > 0. \tag{5.9}$$

*Then the predictive distribution $Q_k$ with predictive density (5.8) is incoherent.*

*Proof.* For each $x \in \mathcal{X}$, the variation distance between $Q_H(\cdot|x)$ and $Q_k(\cdot|x)$ is

$$
\begin{aligned}
&\sup_{B \subseteq R^p} |Q_k(B|x) - Q_H(B|x)| = \\
&\sup_{B \subseteq R^p} \left| \int_B [|L|^{-1} k(L^{-1}z) - |L|^{-1} k_1(L^{-1}z)] dz \right| = \\
&\frac{1}{2} \int |k(z) - k_1(z)| dz.
\end{aligned}
\tag{5.10}
$$

The second equality follows by a simple change of variable and the well known identity involving variation distance (for example, see Billingsley [1], p. 224 for the argument). Since the last expression in (5.10) is positive by assumption, we see that $Q_k(\cdot|x) \neq Q_H(\cdot|x)$ for all $x$. Thus the set $D$ in (5.4) is $\mathcal{X}$ and hence Theorem 5.1 applies. The conclusion follows. $\qquad \square$

As an application of Theorem 5.3, by taking the mean to be zero in a MANOVA model, and applying the results listed as items 2 through 7 in Table 1 of Keyes and Levy [16], one obtains examples of predictive distributions which are incoherent. Predictive distributions 3 through 7 are obtained via a formal Bayes calculation with an improper prior distribution of the form

$$\nu_\beta(d\Sigma) = |\Sigma|^\beta \frac{d\Sigma}{|\Sigma|^{(p+1)/2}}, \tag{5.11}$$

where $\beta$ satisfies $\beta < (n - p + 1)/2$. The restriction on $\beta$ is necessary so the formal Bayes calculation yields a proper posterior for our example. The improper prior (5.11) yields a predictive distribution with a density of the form

$$q_\beta(z|x) = C|L|^{-1} \frac{1}{(1 + (L^{-1}z)'(L^{-1}z))^{(n+1-2\beta)/2}}, \tag{5.12}$$

where $C$ is a constant. Theorem 5.3 implies that all such predictive distributions are incoherent. The details are routine and left to the reader.

**Acknowledgments.** Thanks to David Freedman and a referee for comments on drafts of this paper.